\documentclass[11pt]{amsart}
\usepackage{amsmath}
\usepackage{amscd}
\usepackage{graphics}
\usepackage{latexsym}
\oddsidemargin .in
\theoremstyle{plain}
\newtheorem{Thm}{Theorem}

\newtheorem{Rm}{Remark}

\theoremstyle{remark}

\def\R{{\mathbb R}}

\def\CL{\mathcal L}

\def\IC{\mathbb C}

\def\IP{\bf P}
\def\IR{\mathbb R}
\def\IO{\mathbb O}
\def\IT{\mathbb T}
\def\IZ{\mathbb Z}

\def\s2x{\hbox{$S^2 \times S^2$}}

    \def\sqr#1#2{{\vcenter{\hrule height.#2pt
            \hbox{\vrule width.#2pt height#1pt \kern#1pt
            \vrule width.#2pt}\hrule height.#2pt}}}
    \def\square{\mathchoice\sqr67\sqr67\sqr{2.1}6\sqr{1.5}6}

\begin{document}

\title[]{Mirror Duality in a Joyce Manifold}
\author{Selman Akbulut, Baris Efe, and Sema Salur }
\thanks{First named author is partially supported by NSF grant DMS 0505638}
\address{Department  of Mathematics, Michigan State University, East Lansing, MI, 48824}
\email{akbulut@math.msu.edu }
\address{Department  of Mathematics, Michigan State University, East Lansing, MI, 48824}
\email{efebaris@msu.edu}
\address {Department of Mathematics, University of Rochester, Rochester, NY, 14627 }
\email{salur@math.rochester.edu } 
\date{\today}

\begin{abstract}

Previously the two of the authors defined a notion of dual Calabi-Yau manifolds in a $G_2$ manifold, and described a process to obtain them. Here we apply this process to a compact $G_2$ manifold, constructed by Joyce, and as a result we obtain a pair of Borcea-Voisin Calabi-Yau manifolds, which are known to be mirror duals of each other.

\end{abstract}

\maketitle

\setcounter{section}{0}
\vspace{-0.5in}

\section{Introduction}

Recall that $G_2$ is the simple Lie group which can be identified with
the subgroup 
$$ G_{2}=\{ A \in GL(7,\R) \; | \; A^{*} \varphi_{0} =\varphi_{0}\; \} $$

\noindent where 
$\varphi _{0}
=e^{123}+e^{145}+e^{167}+e^{246}-e^{257}-e^{347}-e^{356}$ with
$e^{ijk}=dx^{i}\wedge dx^{j} \wedge dx^{k}$ (for more information on $G_2$ manifolds the reader can consult \cite{br1}, \cite{br2}, \cite{hl}, \cite{j}, \cite{as2}).
We say a $7$-manifold $M^7$ has a {\it $G_{2}$ structure} if there is  a
3-form $\varphi \in \Omega^{3}(M)$ such that  at each $p\in  M$
the pair $ (T_{p}(M), \varphi (p) )$ is (pointwise) isomorphic to $(T_{0}(
\R^{7}), \varphi_{0})$ (this condition is equivalent to reducing the tangent frame bundle to a $G_2$-bundle).
A manifold with $G_{2}$ structure $(M,\varphi)$  is called a
{\it $G_{2}$  manifold} (integrable $G_2$ structure) if at each point $p\in M$ there is
a chart  $(U,p) \to (\R^{7},0)$ on which $\varphi $ equals to
$\varphi_{0}$ up to second order term, i.e. on the image of the open set $U$ we can write $\varphi (x)=\varphi_{0} + O(|x|^2)$.  
Also the form $\varphi\in \Omega^{3}(M)$ of a manifold with a $G_2$ structure $(M,\varphi)$ induces an orientation, a metric $g=<,>$, and a cross product operation $TM\times TM \mapsto TM$: $(u,v)\mapsto u\times v$ via the relation $\varphi (u,v,w)=<u \times v,w>$. Then the condition that $(M,\varphi )$ be a $G_2$ manifold is equivalent to the condition $d \varphi=d^{*} \varphi =0$.
A 4-dimensional
submanifold $X\subset (M,\varphi )$ is called {\em coassociative } if $\varphi|_X=0$, and a 3-dimensional submanifold $Y\subset M$ is called
{\em associative} if $\varphi|_Y\equiv vol(Y)$; this last condition is
equivalent to the condition $\chi|_Y\equiv 0$,  where $\chi \in \Omega^{3}(M,
TM)$ is the  tangent bundle valued 3-form defined by the
identity: $ \langle \chi (u,v,w) , z \rangle=*\varphi  (u,v,w,z) $. 
\vspace{.05in}



\vspace{1in}

Now recall the construction of \cite{as1}:  If $(M^7, \varphi)$ is a $G_2$ manifold and  $\xi$ is a non-vanishing unit vector field on an open subset of $M$, such that the codimension one distribution $V_{\xi}:= \xi^{\perp}$  is integrable with a leaf $X_{\xi}\subset M$, then on $X_{\xi}$ we can define the following  $2$-form and a complex structure: 
\vspace{.05in}
\begin{equation*}
\omega_{\xi}=\xi \lrcorner \; \varphi \;\;\;\mbox{and} \;\; J_{\xi}(X)=X\times \xi.
\end{equation*}
Also we can define a complex $(3,0)$ form  $\Omega_{\xi} =  \textup{Re}\; \Omega_{\xi} +i \;\textup{Im}\; \Omega_{\xi}$, where
\vspace{.05in}
\begin{equation*}
 \textup{Re}\; \Omega_{\xi} = \varphi|_{V_{\xi}}
\;\;\mbox{and}\;\;\; \textup{Im}\; \Omega_{\xi} = \xi \lrcorner \; *\varphi = \langle \chi, \xi \rangle.
\end{equation*}


{\Thm [\cite{as1}] $X_{\xi}$ induces an almost Calabi-Yau structure  $(X_{\xi},  \omega_{\xi}, \Omega_{\xi}, J_{\xi})$, i.e. $\omega_{\xi}$ is a nondegenerate $2$-form which is co-closed, $J_{\xi} $ is a metric invariant almost complex structure which is  compatible with $\omega_{\xi}$, and $\Omega_{\xi} $ is a non-vanishing $(3,0)$ form with  $Re\;\Omega_{\xi} $ closed. More specifically $\varphi |_{X_{\xi}}= Re\; \Omega_{\xi} $
and  $*\varphi |_{X_{\xi}}= \star \omega_{\xi} $. Furthermore, if
$\CL_{\xi}(\varphi )|_{X_{\xi}}=0$ then $d\omega_{\xi}=0$,  and if
$\CL_{\xi}(*\varphi)|_{X_{\xi}}=0$ then $J_{\xi}$ is integrable; when both of these conditions are satisfied then
$(X_{\xi},\omega_{\xi}, \Omega_{\xi},J_{\xi})$  is a Calabi-Yau
manifold.}

\vspace{.1in}

Note that from  \cite{m}  the condition 
$\CL_{\xi}(*\varphi)|_{X_{\xi } }=0$  (complex geometry of $X_{\xi }$) 
 implies that deforming associative submanifolds of $X_{\xi}$ along $\xi$  in $M$ keeps them associative; and 
$\CL_{\xi '}(\varphi )|_{X_{\xi ' }}=0$ (symplectic geometry of $X_{\xi '}$) implies  that deforming coassociative submanifolds of $X_{\xi '}$ along $\xi '$ in $M$  keeps them coassociative.
Now, if $(M^7, \varphi , \Lambda)$ is a $G_2$ manifold with a non-vanishing  oriented  $2$-frame field $\Lambda =<u,v> $ (every orientable $7$-manifold admits such a $\Lambda$ by \cite{t}),  $\Lambda $ and $\varphi $ gives an associative/coassociative bundle splitting 
\begin{equation*}
T(M)={\bf E}\oplus {\bf V},
\end{equation*} 
where ${\bf E}= <u,v,u\times v>$ and ${\bf V}={\bf E}^{\perp}$ .  
 In \cite{as1}, when $\xi$ and $\xi'$ are two sections lying in ${\bf V}$ and ${\bf E}$  respectively, we called the pairs  $X_{\xi}$ and $X_{\xi'}\;$  {\it dual} to each other (and {\it strong dual} to each other if $\xi $ and $\xi'$ are homotopic through non-vanishing vector fields).   
To emphasize the dependance on the choice of $\Lambda$ (or ${\bf E}$), we will refer  $X_{\xi}$  and $X_{\xi'}$ as dual submanifolds {\it adapted to} $\Lambda$ (or ${\bf E}$).  
 Also it should be noted that, in general the bundle ${\bf E}$ has always a non-vanishing section (because it is a trivial bundle),  whereas {\bf V} may not admit a global non-vanishing section.
 
   \vspace{.05in}

 In \cite{as1}, as an application of Theorem 1, the simplest case of torus $\IT^7$ example was discussed, getting subtori $\IT^{6}_{123567}$ and $\IT^{6}_{234567}$  as dual submanifolds with different complex structures, where subscripts denote the circle factors (note that in \cite{as1}  different circle factors used). In this paper, we apply this theorem to a more nontrivial example of a $G_2$ manifold $(M_{\Gamma},\varphi)$, constructed by Joyce \cite{j} (described in Section 4), and obtain a pair of Borcea-Voisin Calabi-Yau manifolds (\cite{b}, \cite{v}), which are known to be mirror duals of each other; furthermore they each are singular fibrations, in one case with  a complex K3 surface as the regular fiber (7), in the other case with special Lagrangian $3$-torus as the regular fiber (9). It is a curious question whether the mirror dual of a quintic in $\IC\IP^4$ can also be obtained this way?
 
 {\Thm  On the Joyce manifold  $M^7_{\Gamma}$ we can choose an oriented $2$-framed field $ \Lambda$ and corresponding vector fields $\xi$ and $\xi'$ so that the adapted submanifolds $X_{\xi}$ and $X_{\xi'}$ are Borcea-Voisin Calabi-Yau manifolds. In particular they are mirror duals to each other in the sense that $b^{1,1}(X_{\xi})=b^{2,1}(X_{\xi'})$.}

\vspace{.1in}

More specifically, the Joyce manifold used here is  simply connected, and it is constructed by resolving  the quotient of $\IT^7$ by a finite group action $M\to \IT^7/\Gamma $. Hence  if $X_{\xi}$ is a closed submanifold it has to be separating in $M$ (i.e. the complement has two components) and $\xi$ has to vanish in the complement of $X_{\xi}$. In particular, if $\xi$ is non-vanishing on all of $M$ then $X_{\xi}$ can not be a closed submanifold. Here it turns out that for two choices of non-vanishing vector fields $\xi$ and $\xi'$ (on open subset of $M$) we get  $X_{\xi}$ and $X_{\xi'}$ to be Borcea-Voisin submanifolds  (in particular they are closed manifolds), such that  their open dense subsets  $X_{\xi}(t)\subset X_{\xi}$ and $X_{\xi'}(t)\subset X_{\xi'}$  are dual submanifolds of $M$ adapted to an associative distribution ${\bf E}$ of $(M,\varphi)$.
Furthermore, as $t\mapsto 0$, $X_{\xi}(t)$ and  $X_{\xi'}(t)$ converge to $\IT^{6}_{123567}/\Gamma$  and $\IT^{6}_{234567}/\Gamma$ respectively. In short,  Borcea-Voisin manifolds provide natural compactifications of these dual submanifolds inside of the Joyce manifold, and in the limit  they become orbifold quotients of the dual submanifolds of $\IT^7$.  In C.-H. Liu's terminology \cite{l}, these are submanifolds that admit asymptotically coassociative and  associative fibrations (Section 4). 

{\Rm The notion of mirror pairs of Calabi-Yau manifolds originated in the early 90s.  It is part of a duality between two $10$-dimensional string
 theories.  Some years later $11$-dimensional M-theory appeared, and its
 compactification on a $G_2$ manifold is meant to be related to 10-dimensional string theory on a Calabi-Yau. Two different Calabi-Yau manifolds X and X' may have the same SCFT and in this case the complex and symplectic invariants  of X and X' are related. The construction in this paper fits well with these M-theory/string theory dualities. }

\vspace{.05in}

Much of this paper is based on the previous  work of Liu, where he studied various singular fibration structures of Joyce manifolds; here we are revisiting those examples to turn them into applications of \cite{as1}, i.e. finding duals by using Theorem 1. We would like to thank Osvaldo Santillan for suggesting us \cite{l} to construct Borcea-Voisin duals inside Joyce manifolds. For a discussion of this example from physics point of view see \cite{ss}.  Also special thanks to Betul Tanbay for inviting us to IMBM math institute in Istanbul, where this work is done in its inspiring  surroundings.

\section{Revisiting $\IT^7$ example }


Let us start with the $G_2$ structure $(\IT^7, \varphi )$, given by the positive $3$-form: 
\begin{equation}
\varphi=e^{123}+e^{145}+e^{167}+e^{246}-e^{257}-e^{347}-e^{356} 
\end{equation}

\vspace{.1in}
 
By using the notation of \cite{as1}, let us take $T(\IT^7)= {\bf E}\oplus {\bf V}$, where ${\bf E}=\{e_1,e_2,e_3\}$ and ${\bf V}=\{e_4,e_5,e_6,e_7\}$ are the associative/coassociative splitting of the tangent bundle, and let us choose the vector fields $\xi=e_{4}$ and $\xi'=e_{1}$. Then we get $X_{\xi}=\IT^6_{123567}$ and  $X_{\xi'}=\IT^6_{234567}$. The corresponding complex structures $J_{\xi}$, $J_{\xi'}$ on $X_{\xi}$, $X_{\xi'}$   can easily be calculated from (1) by using the definitions $J_{\xi}(X)=X\times \xi $  and $<x\times y,z>=\varphi(x,y,z)$.

  \begin{equation} 
J_{\xi}=
\left(
\begin{array}{ccc}
{\bf e_1}  & \mapsto  & e_5   \\
{\bf e_2} & \mapsto  & e_6  \\
 {\bf e_3} & \mapsto & -e_7
\end{array}
\right)
,
\;\;\;
J_{\xi'}=
\left(
\begin{array}{ccc}
{\bf e_2 } & \mapsto  &- {\bf e_3 }  \\
 e_4 & \mapsto  & -e_5  \\
  e_6 & \mapsto &  -e_7
\end{array}
\right)
\end{equation}

\noindent  These complex structures correspond to the complex coordinates on $X_{\xi} $ and $X_{\xi'}$:

\begin{equation} 
\left(
\begin{array}{cc}
z_1 ={\bf x_1} -i x_5  \\
z_2={\bf x_2}-ix_6  \\
 z_3={\bf x_3} + ix_7
\end{array}
\right)
\;\;\mbox{and}
\;\;\;
\left(
\begin{array}{cc}
w_1={\bf x_2} +i{\bf x_3}  \\
w_2= x_4 +ix_5  \\
w_3=  x_6 +ix_7  
\end{array}
\right)
\end{equation}

\noindent In these expressions the basis of associative bundle ${\bf E}$ is indicated by bold face letters to indicate the differing complex structures on $\mathbb{T}^6$, and also to empasise how these complex structures intermingle these bundles. 

\vspace{.1in}

\section{Borcea-Voisin 3-folds}

\vspace{.1in}

Borcea-Voisin $3$-fold (\cite{b}, \cite{v}) is a Calabi-Yau $6$-manifold $Q^6$, constructed by taking an involution on a $K3$ surface $j: N^4\to N^4$ which acts nontrivially on $H^{2,0}(N)$,  and an involution $i$ of an elliptic curve $E=\IC/\IZ$ inducing  $dz\mapsto -dz$, and by taking the ``small'' resolution $Q$  of the quotient by product involution $k=i\times j$:
$$ Q \to (E\times N)/k$$
Hodge numbers are  $b^{1,1}=11+5a-b$, and 
$b^{2,1}=11+ 5b-a$, where $a$ is the number of components of the fixed point set $Fix(j)$, and $b$ is the total genus of $Fix(j)$. 
In particular,  when $Fix(j)$ is a disjoint union of two $\IT^2$'s, then $Q$ is self-dual (in ``Mirror duality" sense) i.e. 
$b^{1,1}=b^{2,1}=19$ and its  Hodge diamond is given by:

$$
\begin{array}{ccccccc}
&&& 1&&&\\
&&0&&0&&\\
&0&&19&&0&\\
1&&19&&19&&1\\
&0&&19&&0&\\
&&0&&0&&\\
&&& 1&&&\\
\end{array}
$$

Also, the projection  $ E\times N \to E$ induces a singular fibration
$ N \to Q \stackrel{\pi}{\to} S^2 $,
 where $S^2=\IT^2/i$ viewed as an orbifold (the``pillow case'') with four singular points. Over each of these four singular points, the  singular fiber consists of transversally intersecting $4$-manifolds 
$N/j \smile S^2\times T^2  \smile S^2\times T^2$, where the union is taken along the two disjoint $T^2$'s which is $Fix(j)$. This is because the four singular  fibers of the projection $(\IT^2\times N)/k \to S^2$ are $N/j$ , and taking the ``small resolution" of $(\IT^2\times N)/k $ amounts to replacing the  $\IT^2\times cone(\IR\IP^3)$ neighborhood of each fixed point torus with  $\IT^2\times T^{*}\IC\IP^1$ (i.e. the Euler number $-2$ disk bundle over $S^2$).

\vspace{.05in}

  \begin{figure}[ht]  \begin{center}
\includegraphics{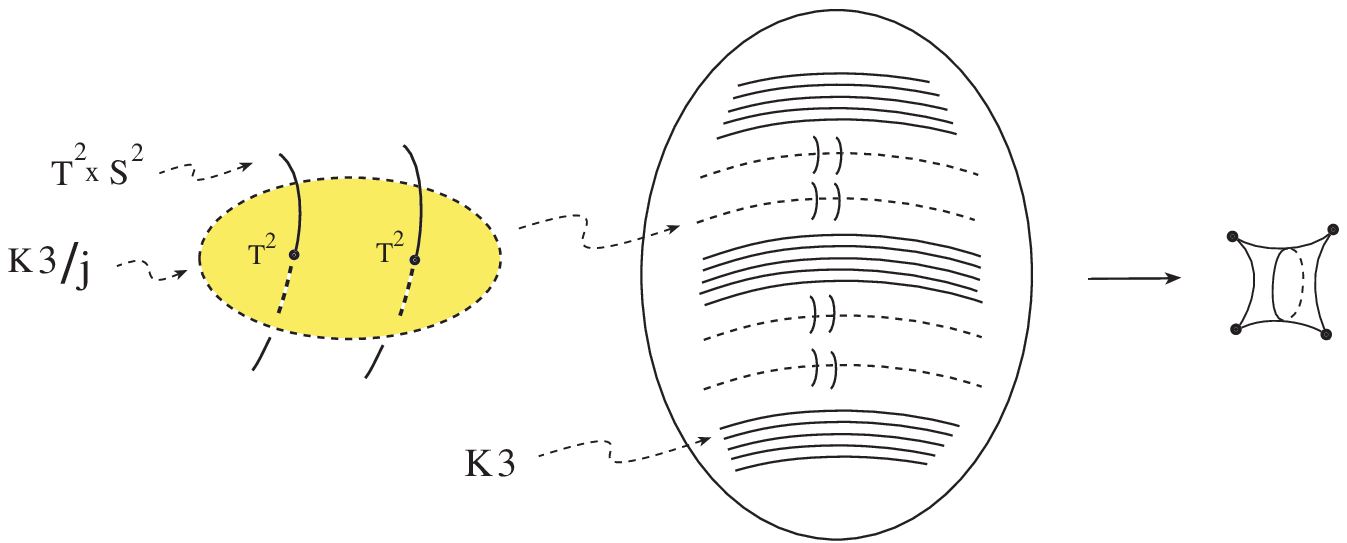}   
\caption{}
\end{center}
   \end{figure}

Also since  $K3/ j \cong \IC\IP^2 \# 9 \bar{\IC\IP^2}$ (e.g.\cite{d}, \cite{dk}), by Van-Kampen theorem $Q$ is simply connected, and hence $Q_{0}:=Q-\{\mbox{point}\}$ is homotopy equivalent to a $4$-complex. Therefore any two  non-vanishing $2$-frame fields on $Q_{0}$  are homotopic to each other, since the obstructions lie in $H^{i}(Q_{0} ,\pi_{i} (V_{7,2}))$ and $V_{7,2}$ is $4$-connected (where $V_{7,2}$ is the Steifel manifold of the oriented $2$-planes in $\IR^7$).



\section{Joyce manifold}


As described in \cite{j}, the group $\Gamma = <\alpha,\beta,\gamma>= \IZ_{2} \oplus \IZ_{2} \oplus \IZ_{2}$ acts on $\IT^7$ by the commuting involutions
$\alpha$, $\beta$ and  $\gamma$   described by the following table:

\begin{table}[h]
\begin{center}
\begin{tabular}{|c||c|c|c||c|c|c|c|} \hline
$  $ &         $x_1$  &  $x_2$ &  $x_3 $ &  $x_4$ & $x_5$ & $x_6$ & $x_7$ \\
\hline \hline
$\alpha$ & $ {\bf x_1}$ & $ x_2$ & $x_3$ & $-x_4$ & $-x_5$& $-x_6$ &$ -x_7$\\
$\beta $ & ${\bf x_1}$  & $ -x_2 $ & $ -x_3 $ & $ x_4$ &$ x_5$ &$-x_6 +\frac{1}{2}$ & $-x_7$\\
\hline
$ \gamma $ &  $- {\bf x_1} $ & $x_2$ & $- x_3$ &  $x_4 $& $-x_5 +\frac{1}{2}$& $x_6 $ & $-x_7 +\frac{1}{2}$\\
\hline
\end{tabular}
\end{center}
\end{table}

\noindent This action was  also studied in \cite{l}. By writing the generators in a different order we observe a certain symmetry with respect to $x_1$ and $x_4$:

\begin{table}[h]
\begin{center}
\begin{tabular}{|c||c|c|c||c|c|c|c|} \hline
$  $ & $x_1$  &  $x_2$ &  $x_3 $ &  $x_4$ & $x_5$ & $x_6$ & $x_7$ \\
\hline \hline
$ \gamma $ &  $-x_1 $ & $x_2$ & $- x_3$ &  ${\bf x_4} $& $-x_5 +\frac{1}{2}$& $x_6 $ & $-x_7 +\frac{1}{2}$\\
$\beta $ & $x_1$  & $ -x_2 $ & $ -x_3 $ & $ {\bf x_4}$ &$ x_5$ &$-x_6 +\frac{1}{2}$ & $-x_7$\\
\hline
$\alpha$ & $ x_1$ & $ x_2$ & $x_3$ & $-{\bf x_4}$ & $-x_5$& $-x_6$ &$ -x_7$\\
\hline
\end{tabular}
\end{center}
\end{table}

\noindent which says that the subgroups $<\alpha, \beta>$ and $<\gamma, \beta>$ induce actions on the $6$-tori $X_{\xi'}$ and $X_{\xi}$, respectively (since they act identity on the remaining coordinate). 
 
 \vspace{.05in}
 
 Each of these involutions fix $16$ disjoint $3$-tori, and all of these $48$ tori (three groups of sixteen) are disjoint from each other. Clearly the fix point sets are:
\begin{eqnarray*}
\mbox{Fix}(\alpha)= \{x_j \in \{0, \frac{1}{2}\}\;\mbox{for}\;\;  j=4,5,6, 7 \} \times \IT^{3}_{123}\\
\mbox{Fix}(\beta)=\{x_j \in \{0, \frac{1}{2}\}\;\mbox{for}\;  j=2,3,7 \;
\mbox{and} \;\;x_6 \in\{\frac{1}{4}, \frac{3}{4}\} \} \times \IT^{3}_{145}\\
\mbox{Fix}(\gamma)=\{ x_j \in \{0, \frac{1}{2}\}\;\mbox{for}\;  j=1,3, \;
\mbox{and} \;\;x_j \in\{\frac{1}{4}, \frac{3}{4}\} \;\mbox{for}\;  j=5,7   \} \times \IT^{3}_{246}
\end{eqnarray*}

\vspace{.05in}

The components of the fix point set of any  one of these involutions are permuted among each other by the other two  in the obvious way;  for example there are four tori $S_1^{\alpha} ,S_2^{\alpha},S_3^{\alpha},S_4^{\alpha} \in \mbox{Fix} (\alpha)$ with
$ \mbox{Fix}(\alpha) =\{ S_j^{\alpha} , \beta (S_{j}^{\alpha}), \gamma (S_{j}^{\alpha}), \beta \gamma(S_{j}^{\alpha}) \;|\; j=1,..,4\}$
Hence the quotient space $\IT^7/\Gamma $ has singular set consisting of  twelve disjoint $3$-tori  $\{S^{\alpha}_{j}, S^{\beta}_{j},S^{\gamma}_{j} \;|\; j=1...4\;\}$, where the neighborhood of each component is described by
\begin{equation} 
 \IT^3\times (\IC^{2}/<-1>)=  \IT^3\times cone(\IR\IP^3).
 \end{equation} 
Then by taking a small resolution along the singular set gives a Joyce manifold $M^7$, i.e. replacing the neighborhood  of components of the singular set of $\IT^7/\Gamma $ by
\begin{equation} 
N=\IT^3\times T^{*}\IC\IP^1
\end{equation} 
gives $M^7$. Here $T^{*}\IC\IP^1$ denotes the disk bundle over $S^2$ with Euler number $-2$. 

\vspace{.1in}

$N$ has a family of torsion free $G_2$ structures depending on a real parameter $t$ 
\begin{equation}
\varphi_{t}=\delta_{1}\wedge \delta_{2}\wedge \delta_{3} +\omega_{1}(t)\wedge\delta_{1}+\omega_{2}(t)\wedge\delta_{2}+ \omega_{3}(t)\wedge\delta_{3}
\end{equation}

\noindent where $\{\delta_{1},\delta_{2},\delta_{3}\}$ are $1$-forms giving the flat structure on $\IT^3$ and $\{\omega_{1}(t),\omega_{2}(t),\omega_{3}(t)\}$ are family of self-dual $2$ forms on $T^{*}\IC\IP^1$. Joyce glues torsion free $G_2$ structure inherited from $\IT^7$ in the complement of the singular set, with this  structure on a small neighborhood $N $ of the divisor, obtains a $G_2$ structure $\varphi_{t}$ which is not torsion free only on a small annular neighborhood of $\partial N$, then he deforms this $G_2$ structure by an exact form $\varphi_{t}\mapsto \varphi_{t} +d\eta_{t}$, making it a torsion free $G_2$ structure \cite{j}.
 Dividing $\IT^7$ by the action of $\Gamma $ and resolving process can be performed in several steps:

\vspace{.15in}
 
$\begin{array}{ccccccccc}
\IT^{7} \hspace{-.2in}&&&&&&& &\\
&\searrow&&&&&& \\
&& \hspace{-.1in}\IT^{7}/<\alpha>&\hspace{-.1in} \leftarrow&  \IT^{3}\times K3 &&&& \\
\hspace{.1in} \downarrow&&\downarrow&&\downarrow&&&& \\
&& \IT^{7}/<\alpha, \beta>&\hspace{-.05in}\leftarrow &  \IT^{1}\times (\IT^{2}\times K3)/<\beta>&\leftarrow &\hspace{-.5in} S^1\times  Q'&&\\
&\swarrow&&&&&\hspace{-.5in} \downarrow&&\\
\IT^{7}/\Gamma \hspace{-.1in} &&&&\leftarrow&&\hspace{-.2in} (S^1\times Q')/<\gamma>&\leftarrow& M\\
\end{array}$

\vspace{.15in} 

\noindent where the horizontal maps are resolutions and vertical maps are quotients by finite group actions. Clearly $\IT^{7}/<\alpha>= \IT^{3}_{123}\times (\IT^{4}_{4567}/\IZ_{2})$ which can be resolved to 
$ \IT^{1}\times \IT^{2}\times K3$ where $K3$ is the Kummer surface,  and the middle diagram commutes. 
The $\beta $ action on $\IT^{7}/<\alpha>$ lifts to an action of $ \IT^{1}\times \IT^{2}\times K3$, where the action $\beta =\mbox{id}\times i\times j$ is identity on the first factor and anti holomorphic in the other two variables, because by (3) the derivative of induced action on $ \IT^{2}\times K3$  is given by: 
$$(w_1, w_2, w_3)\mapsto ( -w_1, w_2, -w_{3}).$$
Hence the resolution of the quotient $(\IT^{2}\times K3)/<\beta>$ is a Borcea-Voisin manifold $Q'$. Then the Joyce manifold $M^7$ is obtained by first folding $S^1\times Q'$ by the $\gamma$ action and getting a singular fibration over the interval $[0,\frac{1}{2}]$, then by resolving the four $3$-tori (the fixed points of the induced action $\gamma$ on the Borcea-Voisin) that lies over  the end points of the interval. So $M^7$ is a fibration on a closed interval $\pi: M^7\to [0,\frac{1}{2}]$ with two singular fibers over the end points $\{0,\frac{1}{2}\}$, and with regular fiber $Q'$.  Each singular fiber is a $6$ manifold with boundary (with singular points)  (e.g. \cite{l}). 

 \begin{figure}[ht]  \begin{center}
\includegraphics{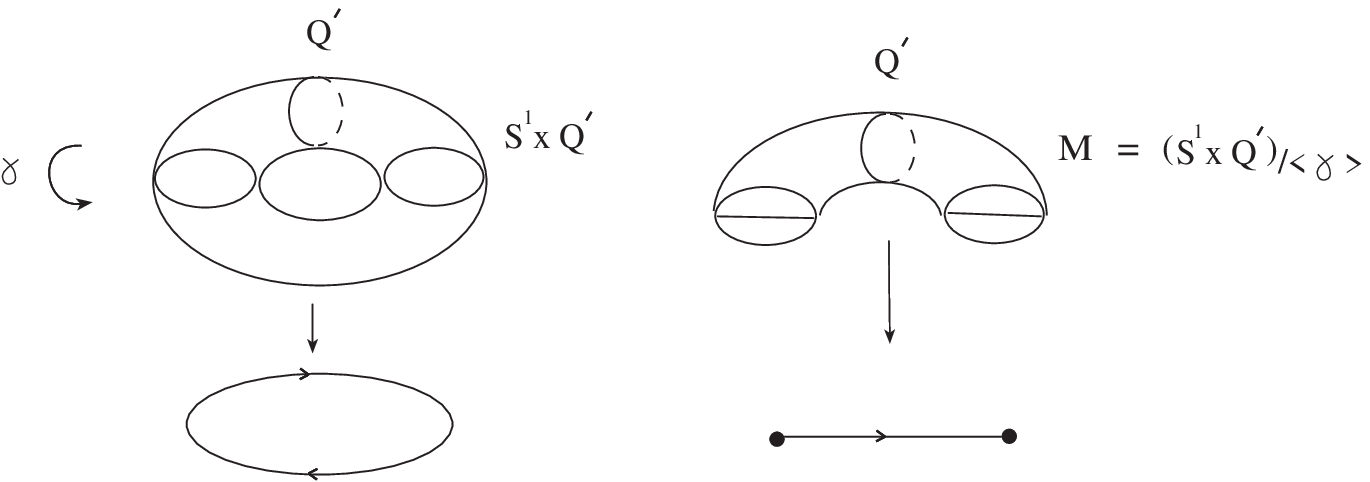}   \caption{}
\end{center}
   \end{figure}
 
  Note that $Q'$ is induced from 
$ \IT^{2}_{23}\times \IT^{4}_{4567}$ inside $\IT^{7}$. More specifically,  $\Gamma$ induces an action of $<\alpha, \beta>$ on a generic slice $ \IT^{6}_{234567}$, fixing $32$  disjoint $2$-tori (two groups of sixteen: the fixed point sets of $\alpha$ and 
$\beta$), and  hence $ \IT^{6}_{234567}/<\alpha, \beta>$ has singular set consisting of sixteen $ \IT^{2}$ (two groups of eight: since $\alpha $ and $\beta$ act on each others fix point sets identifying them in pairs). The neighborhoods of the components of the singular sets are $\IT^2\times 
cone ({\IR\IP^3})$, by resolving them by $\IT^2\times 
T^{*}\IC\IP^1 $ gives $Q'$. It is easy to check that $Fix(j)$ consists of a disjoint union of two $\IT^2$'s, hence $Q'$ is self dual (Section 3). As discussed in Section 3 the projection 
$ \IT^{2}_{23}\times \IT^{4}_{4567} \to  \IT^{2}_{23}$ induces a singular fibration  of $Q'$ over $S^2$, where $S^2$ viewed as the ``pillowcase" orbifold:
\begin{equation}
K3\to Q'\to S^2_{23}.
\end{equation}
Similarly the projection 
$ \IT^{2}_{23}\times \IT^{4}_{4567} \to  \IT^{2}_{4567}$ induces another singular fibration: 
\begin{equation}
\IT^2_{23}\to Q'\to K3/ j.
\end{equation}
with regular fibers $\IT^2$, and singular fibers consisting of union of five $S^2$'s (the pillowcase blown up at four corners) lying over the singular set $Fix(j) \subset K3/j$. Notice, the complex structure (3) $J_{\xi'}$ on $Q'$ makes fibers of the fibrations (7) and (8) complex.

\vspace{.1 in}

 Now in the above discussion, by exchanging the roles of $\alpha$ and $\gamma$ (and $x_1$ and $x_4$) we can do the same analysis (using the second table) and get a fibration on a closed interval $ [0,\frac{1}{2}]$ with two singular fibers over the end points $\{0,\frac{1}{2}\}$, and the regular fibers consisting of a Borcea-Voisin manifold $Q$, which is obtained by resolving the quotient $(\IT^{2}\times K3)/<\beta>$. As above, the Joyce manifold $M$ is now obtained from the fibration $S^1\times Q\to S^1$ by folding it by the $\alpha$ action. From (3) we see that the derivative of the action $\beta $ on $ \IT^{2}\times K3$ is given by: 
$$(z_1, z_2, z_3)\mapsto ( -z_1, z_2, -z_3).$$

\noindent This $Q$ is obtained from  $\IT^{2}_{26}\times\IT^{4}_{1357}$ (also by exchanging the roles of $\gamma$ and $\beta$ here we can get another similar Borcea-Voisin, induced from  $\IT^{2}_{15}\times\IT^{4}_{2367}$). Now the projection 
$ \IT^{2}_{26}\times \IT^{4}_{1357} \to  \IT^{2}_{567}$ induces a singular fibration  of $Q$ over the $3$-sphere, viewed as the orbifold $S^3_{567}= \IT^{3}_{567}/<\beta, \gamma>$,  and $\IT^3$ as the regular fiber.
\begin{equation}
\IT^{3}_{123}\to Q\to S^3_{567}
\end{equation}

To see this first notice that, each of the involutions $<\beta, \gamma >$ fixes four disjoint circles on $ \IT^{3}_{567}$ (and $\beta \gamma$ has no fix points), and also the resulting eight circles are disjoint from each other. The four components of the fix point set of any  one of these involutions are paired among each other, by the other involution (by complex conjugation map), i.e.
\begin{eqnarray*}
\mbox{Fix}(\gamma)=  \{\frac{1}{4}, \frac{3}{4}\}\times S^1_{6} \times  \{\frac{1}{4}, \frac{3}{4}\}\\
\mbox{Fix}(\beta)= S^1_{5}\times  \{\frac{1}{4}, \frac{3}{4}\} \times
 \{0, \frac{1}{2}\}
\end{eqnarray*}

\noindent Hence the singular set of $S^3_{567}$ consists of four disjoint circles (because $\alpha$ and $\beta$ identified pair of circles from each other's fixed point sets). Now we can see $S^3_{567}$ as a $3$-sphere, by first taking the quotient of $\IT^{3}_{567}/<\beta > =  \IT^1_{5} \times S^2_{67} $ (circle cross the ``pillowcase''), then by folding by the action $\gamma$.
 
\begin{figure}[ht]  \begin{center}
\includegraphics{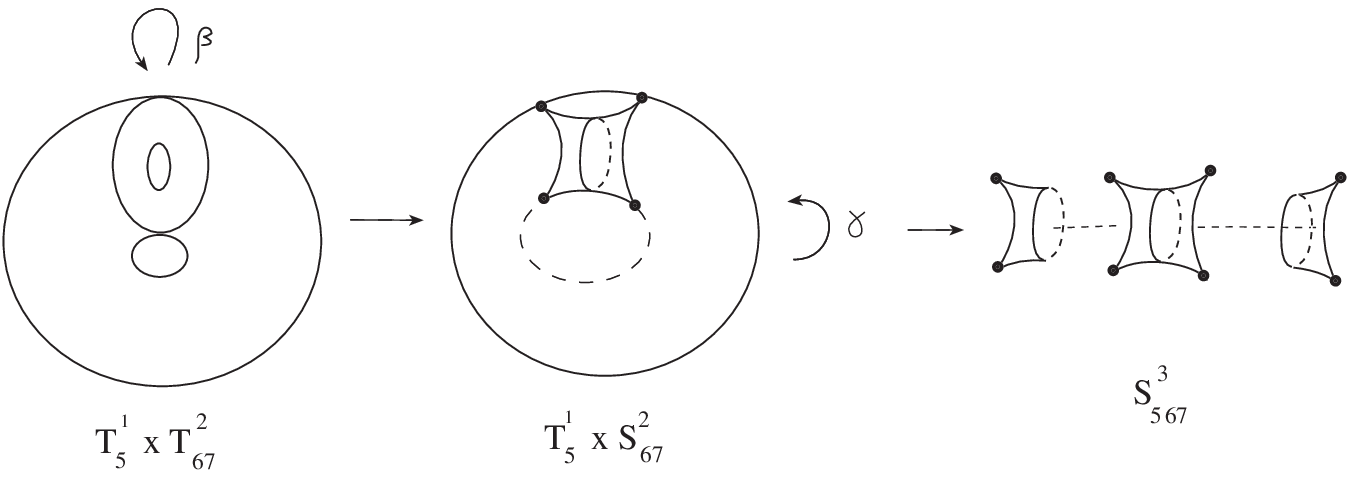}   \caption{}
\end{center}
   \end{figure}
   
\noindent  As noted in \cite{l}, from this one can see that the singular set of the orbifold $S^3_{567}$  is the four component link $L\subset S^3$ shown in Figure 4.  As an independent  check, by using the methods of \cite{ak} the reader can draw the handlebody picture of the $4$-fold ${\bf Z}_{2}\oplus {\bf Z}_{2}$ - branched covering of $S^3$ branched along $L$ (with branching index $2$ on each component of $L$) and get $T^3$.   
   \begin{figure}[ht]  \begin{center}
\includegraphics{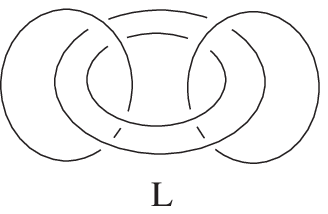}   
\caption{}
\end{center}
   \end{figure}

  Also, it  follows from the identifications that the singular fiber of the fibration \begin{equation}
\IT^6_{123567}/<\beta, \gamma> \to S^3_{567}
\end{equation}
 is $S^2\times S^1$, that is the points lying over the link $L$,  this is because
$$\IT^3_{123}/<\gamma>=\IT^3_{123}/(x_1,x_2,x_3)\sim (-x_1,x_2,-x_3) =S^{2}_{13}\times S^1.$$
Therefore the  fibration (9) is obtained by resolving $\IT^3_{123567}/<\beta, \gamma>$ along the sixteen singular $T^2$'s = $\{$four corners of $S^{2}_{13} \} \times S^1 \times L$, i.e. by  replacing the singular neighborhood $cone(\IR\IP^3)\times \IT^2$ by $ T^{*}{\IC\IP^1}\times \IT^2 $. Also note that the complex structure (3) $J_{\xi}$ on $Q$ makes the regular fiber of the fibration (9) special Lagrangian.
 
 \vspace{.1in}

  Now recall that the Joyce manifold $M$ is fibered over the interval in two different ways (Figure 2). The obvious vector field on the interval with zeros at the end points lifts to transverse vector fields $\xi'$ and $\xi$ to the generic fibers of $X_{\xi'}$ and $X_{\xi}$  of these two fibrations (they were called $Q'$ and $Q$ above). Clearly $\xi'$ and $\xi$  descend from the vector fields $e_1$ and $e_4$ of $\IT^7$. 
  
 \vspace{.05 in} 
 
  From the fibrations  (8) and (9), and from (6),   we see that the subbundle  generated by   $\{e_1,e_2,e_3\}\subset \IT^7$ descends as an associative subbundle  ${\bf E}$ on  a neighborhood $\IO_{t}$  of the regular fibers of $Q$ and $Q'$ inside $M$. Note that in the case of (8) $e_1$ is the transverse direction to $Q'$ in $M$. Say $X_{\xi}(t)=\IO_{t}\cap X_{\xi}$ and $X_{\xi'}(t):=\IO_{t}\cap X_{\xi'}$  are the complements of small radius $t$-tubes around the singular fibers. It is clear that, by the terminology of Section 1,  $X_{\xi}(t)$ and $X_{\xi'}(t)$ are dual submanifolds adapted to ${\bf E}$, and as $t\to 0$ they converge to the orbifold quotients $\IT^{6}_{123567}/\Gamma$  and $\IT^{6}_{234567}/\Gamma$. 


{\Rm It is an interesting question whether the $G_2$ structure $\varphi$ on $M$ given by Joyce's theorem can be deformed so that either one of  the conditions  $\CL_{\xi}(\varphi )|_{X_{\xi}}=0$ or $\CL_{\xi}(*\varphi )|_{X_{\xi}}=0$ hold. If this is the case then  we can conclude that the abstract Borcea-Voisins $X_{\xi}$ are indeed  exactly  symplectic or complex Calabi-Yau submanifolds by the induced structure given by the Theorem 1 (not just asymptotically). } 

{\Rm If $\lambda $ is a $2$-frame in the vertical tangent bundle of the fibrations (8) and (9) (i.e. tangent vectors of the regular fibers), by the remarks at the end of the Section 3, we may assume that $X_{\xi}(t)$ and $X_{\xi'}(t)$ are duals adapted to a global non-vanishing $2$-frame field $\Lambda$ on $M$. This is bacause we can first choose a non-vanishing $2$-frame field $\Lambda'$ on $M$ by using \cite{t}, then over $X_{\xi}(t)\cup X_{\xi'}(t)$ we can homotope the $2$-frame $\lambda $ to $\Lambda' $, then by the homotopy extension property we can extend  $\lambda $  over $M$ as a non-vanishing $2$-frame $\Lambda$.}

{\Rm Some of the $Spin(7)$ manifolds constructed in \cite{j} from the quotients of $\IT^8$ (e.g. Example 14.2), also fiber over an interval in two different ways with dual ``almost $G_2$ manifolds''  (in the sense of \cite{as1}) as generic fibers. So we might hope to obtain more interesting ``dual" Calabi-Yau's living in two different  $G_2$'s which are themselves ``dual'' in a $Spin(7)$; in the case of Example 14.2 of \cite{j} again we get a pair of Borcea-Voisin's. }


\end{document}